\documentclass[10pt,twoside]{siamltex}
\usepackage{amsfonts,epsfig}
\usepackage{amsmath,amssymb}
\usepackage{amsfonts}
\newcommand{\ef}{\mathbb{F}}
\newcommand{\enn}{\mathbb{N}}

\newcommand{\Hom}{\mathbf{Hom}}
\newcommand{\Homf}{\mathbf{Hom}_{\ef}}

\begin{document}

\bibliographystyle{plain}
\title{Duality of discrete topological vector spaces}
% Leave blank; editors will write the exact dates above

\author{
Ramamonjy ANDRIAMIFIDISOA.\thanks{D\'{e}partement de Math\'{e}matiques et Informatique,
Facult\'{e} des Sciences, BP 906, Universit\'{e} d'Antananarivo, MADAGASCAR
(ramamonjy.andriamifidisoa@univ-antananarivo.mg)}
% Remember to put \and between any two authors
}
% Note that \footnotemark[3]} is used for the third author
% because of the same affiliation for the second and third authors.
% If the same affiliation is to be used for the first and second authors,
% \footnotemark[2] should be used instead of \thanks{} for the second author.

% Authors and running title to go on top of each page
\pagestyle{myheadings}
\markboth{Ramamonjy ANDRIAMIFIDISOA}{Duality of discrete topological vector spaces}
\maketitle

\begin{abstract}For a field $\ef$, the discrete topological vector spaces over $\ \ef$ are essentially of  the form $\ef^{\alpha}$ where $\alpha$ is an ordinal. With additional appropriate properties, they are isomorphic to $\ef^{(\beta)}$ where $\beta$ is again an ordinal. Finally, the categories of the vector spaces of the the first and the second type are equivalent. \end{abstract}

\begin{keywords}
Discrete topological vector space, ordinal, inverse limits, adjoint of a linear mapping, duality of categories
\end{keywords}
\begin{AMS}15A03, 15A04, 18A30, 18B30
\end{AMS}

%%%%%%%%%%%%%%%%%%%%%%%%%%%%%%%%%%%%%%%%%%%%%%%%%%%%%%%%%%%%%
\section{Introduction}Let $\ef$ be a field and $\enn$ the set of the natural integers. We write $\ef^{\enn}$  for the set of all sequences $y :\enn\longrightarrow\ef, i\longmapsto y_i$ and $\ef^{(\enn)}$ for the subset of $\ef^{\enn}$ consisting of the sequences $y$ with finite support. These are vector spaces over $\ef$. The field $\ef$ is endowed with the \textit{discrete topology} and  becomes a topological vector space. With the product topology, $\ef^{\enn} $ becomes a topological vector space too, whose $0$-basis consists  of sets $(V_n)_{n\in\enn^{*}}$ where
for all $n\in\enn^*$,
\begin{equation*}
V_n = \{y\in\ef^{\enn}\;\vert\;y_i = 0\quad\text{for all}\quad i<n\}.
\end{equation*}
It follows  immediately that
\begin{equation*}
V_1  \supset V_2\supset\ldots \supset V_n \supset V_{n+1}\supset\ldots\\
\end{equation*}
Let $\mathfrak{G}$ be the category of  the discrete $\ef$-vector spaces : an object of  $\mathfrak{G}$ is a vector space with a denumerable algebraic basis. If $U$ and $V$ are two objects of $\mathfrak{G}$, a morphism $f :U\longrightarrow V$ is  just a linear map from $U$ to $V$. We write $\Hom_{\mathfrak{G}}(U,V)$ for the set of morphisms from $U$ to $V$.\\
The category $\mathcal{T}\mathfrak{G}$  is the subcategory of $\mathfrak{G}$ whose objets are the topological (discrete) vector spaces $(W,\mathcal{T}_W)$ verifying the following properties :\\
$(1)$\;\;A basis   of the filter of neighbourhoods of $0$ consists of the lattice $\mathcal{C}$ of the subspaces $C$ of
 $W$ of finite codimension,\\
 $(2)$\;\;The topology $\mathcal{T}_W$ is Haussdorf, i.e
 \begin{equation*}
\bigcap_{C\in\mathcal{C}} = \{0\},
\end{equation*}
$(3)$\;\;Le lattice $\mathcal{C}$ has a denumerable subset $\mathcal{B}=\{V_n\;\vert\;n\in\beta\}$ (where $\beta$ is an ordinal) which verifies
\begin{equation*}
W = V_0\supset V_1\supset\ldots\supset V_n\supset\ldots\quad\text{with}\quad\dim(W/V_n) = n.
\end{equation*}
$(4)$\;\;The topology $\mathcal{T}_W$ is complete, i.e
\begin{equation*}
W = \underleftarrow{\lim}(W/V_n).
\end{equation*}
For every two objects $U$ and $V$ of $\mathcal{T}\mathfrak{G}$, a morphism $f :U\longrightarrow V$ is a continuous linear mapping. We write $\Hom_{\mathcal{T}\mathfrak{G}}(U,V)$ for the set of morphisms from $U$ to $V$.\\
The functor $\Homf(-,\ef)$ is defined by
\begin{align}\label{homf}
\begin{split}
\Homf(-,\ef) :\mathfrak{G}&\longrightarrow\mathcal{T}\mathfrak{G}\\
U&\longmapsto\Homf(U,f)\\
 \Hom_\mathfrak{G}(U,V)&\longmapsto\left\{
                                            \begin{array}{ll}
                                              \Homf(f,\ef)\in\Hom_{\mathcal{T}\mathfrak{G}}(\Homf(V,\ef),\Homf(U,\ef))\\ \\
                                             \Homf(f,\ef)(y) =y\circ f\quad\text{for all}\quad y\in\Homf(V,f).
                                            \end{array}
                                          \right.
\end{split}
\end{align}
Our goal is to show that this functor induces a duality (\cite{andriam-1}, \cite{lane}, \cite{oberst}) between the categories $\mathfrak{G}$ and $\mathcal{T}\mathfrak{G}$.  For this purpose, parting from $\ef^{\enn}$, we introduce the topological vector spaces $\ef^{\alpha}$ and $\ef^{(\alpha)}$ where $\alpha$ is an ordinal. Using inverse limits of sets and categories (\cite{bou}),  we characterise $\mathfrak{G}$ and $\mathcal{T}\mathfrak{G}$ in section \ref{dtvs}. The main result is stated in theorem \ref{main}.
\section{The topological vector space of sequences over a field}\label{dtvs}We begin with a lemma, which characterizes of $\ef^{\enn}$ :

\lemma\label{lem1} The following properties hold :\\
$(1)$\;\;$\ef^{\enn}/V_n\simeq\ef^n$ for all $n\in\enn^*$.\\
$(2)$\;\;For every $m,n\in\enn^*$ such that $m\leqslant n$, let
\begin{align}\label{gmn}
\begin{split}
g_{m,n}:\ef^{\enn}/V_n&\longrightarrow\ef^{\enn}/V_m\\
\bar x^n&\longrightarrow\bar x^m
\end{split}
\end{align}
where for all $x\in\ef^{\enn}$,  $\bar x^k$ is the class of $x$ modulo $V_k$. Then
\begin{equation}\label{liminv}
 \underleftarrow{\lim}(\ef^{\enn}/V_n)\simeq\ef^{\enn}
\end{equation}
where $\underleftarrow{\lim}(\ef^{\enn}/V_n)$ is the \textbf{inverse limit} of the sequence $(\ef^{\enn}/V_n)_{n\in\enn^*}$ with respect to the $g_{m,n}$
\endlemma
\begin{proof} $(1)$\;\;For all $n\in\enn^*$, we have
\begin{equation*}
\ef^{\enn} = V_n\oplus U_n
\end{equation*}
where
\begin{align*}
\begin{split}
U_n= \{y\in\ef^{\enn}\;&\vert\;y_i = 0 \quad\text{for all}\quad i\geqslant n\}\\
 &=\{(y_1,y_2,\ldots,y_n,0,\ldots,0,\ldots)\;\vert\;y_i\in\ef\quad\text{for all}\quad i= 1,\ldots,n\}.
\end{split}
\end{align*}
Thus $U_n\simeq\ef^n$ and $\ef^{\enn}/V_n\simeq U_n\simeq\ef^n$. \\
$(2)$\;\;Given  $m,n\in\enn^*$ with $m\leqslant n$, consider the mapping

\begin{align*}
\begin{split}
f_{m,n} : \ef^n&\longrightarrow\ef^m\\
(x_1,\ldots,x_m,x_{m+1},\ldots,x_n)&\longmapsto(x_1,x_2,\ldots,x_m).
\end{split}
\end{align*}
By definition, the \textit{inverse limit} (\cite{bou}) of the sequence $(\ef^n)_{n\in\enn^*}$ with respect to these mapping is
\begin{align}\label{invlim}
\begin{split}
&\underleftarrow{\lim}(\ef^n)= \{(x_1,\ldots,x_{n+1},\ldots)\in\prod_{n\in\enn^*}\ef^n\;\vert\;f_{m,n}(x_n) =  x_m\;\;\forall m,n\in\enn^*\;\;,m\leqslant n\}\\
& = \{(x_1,(x_1,x_2),\ldots,(x_1,\ldots,x_{n-1}),(x_1,\ldots,x_n),(x_1,\ldots,x_{n+1}),\ldots)\}\subset\prod_{n\in\enn^*}\ef^n,
 \end{split}
\end{align}
so that, on the one hand,
\begin{equation*}
\underleftarrow{\lim}(\ef^n) \simeq\ef^{\enn},
\end{equation*}
and on the other hand, by the isomorphism
\begin{align*}
\begin{split}
\ef^{\enn}/V_n&\longrightarrow\ef^n\\
\bar x^n&\longrightarrow (x_1,\ldots,x_n),
\end{split}
\end{align*}
for any $n\in\enn^*$, we get, by taking the inverse limit of the sequence $(\ef^{\enn}/V_n)_{n\in\enn}$ with respect to the mappings $f_{m,n}$ :
\begin{equation*}
\underleftarrow{\lim}(\ef^n/V_n) =\underleftarrow{\lim}(\ef^n)\simeq\ef^{\enn}.
\end{equation*}
\end{proof}

Now we state the main result of the section
\proposition\label{ob} $(1)$\;\;Any object of $\mathcal{T}\mathfrak{G}$ (resp. of $\mathfrak{G}$) is isomorphic to some $\ef^\alpha$ (resp.  $\ef^{(\alpha)}$) where $\alpha$ is an ordinal.\\
$(2)$\;\;Any morphism $f :\ef^{(\beta)}\longrightarrow\ef^{(\alpha)}$ is of the form
\begin{equation*}\label{eq1}
x\longmapsto x\cdot\mathcal{F}
\end{equation*}
where $\mathcal{F}\in\ef^{\beta,(\alpha)}$.\\
$(3)$\;\;Any morphism $g :\ef^\alpha\longrightarrow\ef^\beta$ is of the form
\begin{equation*}\label{eq2}
y\longmapsto\mathcal{G}\cdot y
\end{equation*}
where $\mathcal{G}\in\ef^{\beta,(\alpha)}$.
\endproposition
\proof (1)\;\;For the category $\mathcal{T}\mathfrak{G}$ :
 for $W\in\mathcal{T}\mathfrak{G}\setminus\{0\}$, if $\dim W= n<+\infty$, then $W\simeq\ef^n$. So, assume that $\dim W=+\infty$; it is known that there are $n\in\enn^*$ and a subspace $U_n$ of $W$ such that
\begin{equation*}
W = V_n\oplus U_n\;\;\text{with}\;\; W/V_n\simeq U_n,\;\;\dim U_n = n\;\;\text{and}\;\; U_m\subset U_n\;\;\text{for all}\;\; m\leqslant n.
\end{equation*}
Now, consider the category whose objets are
respectively the sets $(W/V_n), U_n$ et $\ef^n$ endowed, for all $m\leqslant n$, with the morphisms $f_{m,n}$ and $g_{m,n}$ as in the lemma (\ref{lem1}) and its proof and the projections $h_{m,n} : U_n\longrightarrow U_m$ of $U_n$ on $U_m$. By taking the inverse limits with respect to these categories and using  the lemma (\ref{lem1}), we obtain
\begin{equation*}
W = \underleftarrow{\lim}(W/V_n) \simeq\underleftarrow{\lim}(U_n)\simeq\ef^{\enn}.
\end{equation*}
For the category $\mathfrak{G}$ :
If $P\in\mathfrak{G}$ and $\dim P= n<+\infty$, then $P\simeq\ef^n=\ef^{(n)}$. So, assume that $\dim P=+\infty$. Le $(b_n)_{n\in\enn}$ a basis of $P$; any element $v\in P$ can be uniquely  expressed as
\begin{equation*}
v = \sum_{n\in\enn}\lambda_n v_n
\end{equation*}
where the sequence  $(\lambda_n)_{n\in\enn}$  is with finite support, i.e   an element of $\ef^{(\enn)}$ . Thus $P\simeq\ef^{(\enn)}$.\\
$(2)$\;\;Any morphism $f :\ef^{(\beta)}\longrightarrow\ef^{(\alpha)}$ is defined by its matrix with respect to the canonical bases $\delta=(\delta_i)_{i\in\beta}$ of $\ef^{(\beta)}$ and $\delta'=(\delta'_j)_{j\in\alpha}$ of $\ef^{(\alpha)}$. This matrix is of the form
\begin{equation*}
\mathcal{F} = (d_{ji})_{j\in\beta,i\in\alpha}\quad\text{with}\quad d_{ji}\in\ef\quad\text{for all}\quad (j,i)\in\beta\times\alpha.
\end{equation*}
Since for any $j\in\beta$, we have
de definite sum\begin{equation*}
\mathcal{F}(\delta_j) = \sum_{i\in\alpha}d_{ji}\delta'_i,
\end{equation*}
the sequence  $(\delta_{ji})_{i\in\alpha}$ is necessarily with finite support. Thus $\mathcal{F}\in\ef^{\beta\times(\alpha)}$ and conversely, such a matrix defines a morphism $f:\ef^{(\beta)}\longrightarrow\ef^{(\alpha)}$ by the equation \eqref{eq1}.\\
$(3)$\;\;Let $\mathcal{G}$ be defined by the equation \eqref{eq2}. For any $y\in\ef^{\alpha}$, we have, for all $j\in\beta$
\begin{equation*}
(\mathcal{G}\cdot y)_j = \sum_{i\in\alpha}\mathcal{G}_{ji}y_i
\end{equation*}
and since this sum is definite, the sequence $(\mathcal{G}_{ji})_{i\in\alpha}$ is necessarily with finite support, i.e $\mathcal{G}\in\ef^{\beta,(\alpha)}$, and  such a matrix defines  a morphism $g :\ef^{\alpha}\longrightarrow\ef^{\beta}$ by the equation \eqref{eq2}.  This morphism is continuous in $0$, therefore continuous on $\ef^\alpha$.\\
Conversely, let $g :\ef^{\alpha}\longrightarrow\ef^{\beta}$ be a continuous linear mapping. If $\alpha\in\enn$, then the existence of $g$ in the equation \eqref{eq2} is trivial, so assume that $\alpha=\enn$. For any $y\in\ef^{\enn}$, we have $y=\lim y^{(n)}$ for the topology of $\ef^{\enn}$ where
\begin{equation*}
y^{(n)} = (y_1,y_2,\ldots,y_n,\ldots,0,\ldots)\in\ef^{(\enn)}.
\end{equation*}
Using the continuity of $g$, we have
\begin{equation*}
g(y) = \lim g(y^{(n)}).
\end{equation*}
By the injections $\ef^{(n)}\hookrightarrow\ef^{(\enn)}\hookrightarrow\ef^{\enn}$,
we may see $\ef^{(n)}$ as a subspace of $\ef^{\enn}$. Let $g(n)$ be the restriction of $g$ to $\ef^{(n)}$; it is defined by a matrix $\mathcal{G}^{(n)}=(\mathcal{G}^{(n)})_{ji}\in\ef^{\beta, n}$, i.e.
\begin{equation*}
g^{(n)}(y^n)  = \mathcal{G}^{(n)}\cdot y^{(n)}.
\end{equation*}
Now, for any $i\in\enn$ such that $i\leqslant n$,we have
\begin{equation*}
g(\delta_i) = g^{(n)}(\delta_i)\cdot \mathcal{G}^{(n)}(\delta_i)
\end{equation*}
and for any $j\in\beta$,
\begin{equation*}
g(\delta_i) _j = \sum_{\nu\in\enn}\mathcal{G}^{(n)}_{j\nu}\cdot(\delta_i)_\nu = \mathcal{G}^{(n)}_{ji}.
\end{equation*}
It follows that, for any $j\in\beta$,
\begin{equation*}
g^{(n)}(y^{(n)}_j) = \sum_{i\in\enn}\mathcal{G}_{ji}^{(n)}\cdot y^{(n)} = \sum_{i=1}^n(g(\delta_i))_j\cdot y_i.
\end{equation*}
Hence the sequence $(\sum_{i=1}^n(g(\delta_i))_jy_i)_{n\in\enn}$ converges to $(g(y))_j$ in $\ef$. Since the topology of $\ef$ is discrete,this implies that the sequence   $(\sum_{i=1}^n(g(\delta_i))_jy_i)_{n\in\enn}$ is stationary. Thus the sequence $(g(\delta_i)_j)_{j\in\enn}$ is necessarily with finite support. Therefore, there exists $N\in\enn$ such that
\begin{equation*}
g(y)_j = \lim\sum_{i=1}^ng((\delta_i))_jy_i = \sum_{i=1}^{N_j}g(\delta_{i})y_i.
\end{equation*}
Now, let $\mathcal{G}\in\ef^{\beta,\alpha}$ be the matrix defined by $\mathcal{G}_{ji}=g(\delta_i)_j$. Then $\mathcal{G}\in\ef^{\beta,(\alpha)}$ and
\begin{equation*}
g(y) = \mathcal{G}\cdot y\quad\text{for all }\quad y\in\ef^{\enn}.
\end{equation*}
\endproof

\section{Duality}

In the equation (\ref{homf}), we must show that the functor $\Homf(f,\ef)$ is an element of $\Hom_{\mathcal{T}\mathfrak{G}}(\Homf(V,\ef),\Homf(U,\ef)$. We may take $U=\ef^{\alpha}$ and $V=\ef^{\beta}$ with $\alpha,\beta\in\enn$. For any $y\in\ef^{\beta}$, the matrix of the linear mapping $f :\ef^{\beta}\longrightarrow\ef$ is $\mathcal{H}\in\ef^{\beta,1}$ with
\begin{equation*}\label{hhom}
y(\delta_i) = y_i =\delta_i\cdot\mathcal{H} =\sum_{j\in\beta}(\delta_i)_j\cdot\mathcal{H}  = \mathcal{H}_{i1}
\end{equation*}
for any $i\in\beta$. If $\mathcal{F}$ is the matrix of  $f :\ef^{(\beta)}\longrightarrow\ef^{(\alpha)}$, then the matrix of the linear mapping $y\circ f$  is $\mathcal{H}\circ \mathcal{F}$.
Thus it is a continuous linear mapping with respect to the topology $\mathcal{T}\mathfrak{G}$.\\

The functor $\Homf(-\ef)$ \textit{is exact}, i.e. transforms an exact sequence of vector spaces of $\mathfrak{G}$
\begin{equation*}
0\longrightarrow U\stackrel{f}{\longrightarrow}V
\end{equation*}
 to the exact sequence of vector spaces of $\mathcal{T}\mathfrak{G}$
 \begin{equation*}
0\longleftarrow\Homf(U,\ef)\stackrel{\Homf(f,\ef)}{\longleftarrow}\Homf(V,\ef).
\end{equation*}
This is because any subspace of a vector space admits a supplementary vector space and a linear mapping is defined by its restrictions to all its supplementary subspaces.
\proposition\label{adjoint}Let $\alpha,\beta\in\enn$ and a linear mapping
\begin{align}\label{Matr}
\begin{split}
f :\ef^{(\beta)}&\longrightarrow\ef^{(\alpha)}\\
x&\longmapsto x\cdot \mathcal{F}
\end{split}
\end{align}
where $\mathcal{F}\in\ef^{\beta,1}$. Then $\mathcal{F}$ is also the matrix of $\Homf(f,\ef)$, i.e.
\begin{align*}
\begin{split}
\Homf(f,\ef) :\ef^{\alpha}&\longrightarrow\ef^{\beta}\\
y&\longmapsto\mathcal{F}\cdot y.
\end{split}
\end{align*}
\endproposition
\proof Applying the  functor $\Homf(-,\ef)$ to the equation \eqref{Matr}, we obtain
\begin{align*}
\begin{split}
\Homf(f,\ef) :\ef^{\alpha}&\longrightarrow\ef^{\beta}\\
y&\longmapsto y\circ f,
\end{split}
\end{align*}
and knowing that there is a matrix $\mathcal{G}\in\ef^{\beta,(\alpha)}$ such that $\Homf(f,\ef)(y)=\mathcal{G}\cdot y$ for any $y\in\ef^{\beta}$, we get
\begin{equation}\label{gyf}
\mathcal{G}\cdot y = y\circ f.
\end{equation}
But, by the above equation, the matrix of the linear mapping $y :\ef^{\beta}\longrightarrow\ef$
is $\mathcal{H}\in\ef^{\beta,1}$ with $\mathcal{H}_{ji}=y_j$ for any $j\in\beta$. Therefore, we can write the equation \eqref{gyf} under the following matrix form
\begin{equation*}
\mathcal{G}\cdot\mathcal{H} = \mathcal{F}\cdot\mathcal{H}.
\end{equation*}
Since this is true for any matrix $\mathcal{H}\in\ef^{\beta,1}$, we finally conclude that $\mathcal{G} = \mathcal{F}$.\\ Let $\alpha,\beta\in\enn$ and \begin{align*}
\begin{split}
f:\ef^{(\beta)}&\longrightarrow\ef^{(\alpha)}\\
x&\longmapsto x\cdot\mathcal{F},
\end{split}
\end{align*}
be a linear mapping, where $\mathcal{F}\in\ef^{\beta,1}$. Let
\begin{align*}
\begin{split}
f':\ef^{\beta}&\longrightarrow\ef^{\alpha}\\
y&\longmapsto\mathcal{F}'\cdot y
\end{split}
\end{align*}
be another linear mapping. We say that $f$ and $f'$ are \textbf{adjoints} if $\mathcal{F}'=\mathcal{F}$.
We also say that $f$ (resp. $f'$) is the adjoint of $f$' (resp. $f$).  The adjoint always exists and is unique because it has the same matrix as the given linear mapping : By proposition \ref{adjoint}, for any  $U,V\in\mathfrak{G}$, the adjoint of an element $f\in\Hom_\mathfrak{G}(U,V)$  is the morphism $\Homf(f,\ef)$ of $\Hom_{\mathcal{T}\mathfrak{G}}(\Homf(U,\ef),\Homf(V,\ef))$.
\theorem\label{main} The functor $\Homf(-,\ef)$ induces a duality between the categories $\mathfrak{G}$ and $\mathcal{T}\mathfrak{G}$.
\endtheorem
\proof By theorem (\ref{ob}), we know that any object of $\mathcal{T}\mathfrak{G}$ is isomorphic to an object of $\mathfrak{G}$. Therefor, is suffices to show that the functor $\Homf(-,\ef)$ is \textit{faithful} and \textit{full}, i.e., for all $U,V\in\mathfrak{G}$, the mapping
\begin{align*}
\begin{split}
\Hom_\mathfrak{G}(U,V)&\longrightarrow\Hom_{\mathcal{T}\mathfrak{G}}(\Homf(V,\ef),\Homf(U,\ef))\\
f&\longmapsto\Homf(f,\ef)
\end{split}
\end{align*}
is injective and surjective (then bijective). We can, without loosing generality, take $U=\ef^{(\alpha)}$ and $V=\ef^{(\beta)}$ where $\alpha,\beta\in\enn$.
\\
Surjectivity : we know that an element of $\Hom_{\mathcal{T}\mathfrak{G}}(\Homf(V,\ef),\Homf(U,\ef))$ is of the following form
\begin{align*}
\begin{split}
F :\ef^\alpha&\longrightarrow\ef^\beta\\
y&\longmapsto\mathcal{G}\cdot y
\end{split}
\end{align*}
where $\mathcal{G}\in\ef^{\beta,(\alpha)}$. Then the mapping
\begin{align*}
\begin{split}
f :\ef^{(\beta)}&\longrightarrow\ef^{(\alpha)}\\
x&\longmapsto x\cdot\mathcal{G}
\end{split}
\end{align*}
is well defined and  is an element of $\Hom_\mathfrak{G}(U,V) $ which  verifies $\Homf(f,\ef)=F$.\\
Injectivity : an element $f\in\Hom_\mathfrak{G}(U,V)$ is of the form
\begin{align*}
\begin{split}
f :\ef^{(\beta)}&\longrightarrow\ef^{(\alpha)}\\
x&\longmapsto x\cdot\mathcal{F}
\end{split}
\end{align*}
with $\mathcal{F}=(\mathcal{F}_{ji})_{ji}\in\ef^{\beta,(\alpha)}$ and an element $\varphi\in\Homf(V,\ef)$ is of the form
\begin{align*}
\begin{split}
\varphi\ :\ef^{(\alpha)}&\longrightarrow\ef\\
y&\longmapsto \mathcal{G}\cdot y
\end{split}
\end{align*}
where $\mathcal{G}=(\mathcal{G}_{ij})_i\in\ef^{\alpha,1}$. Then we have
\begin{align*}
\begin{split}
\varphi(f(x)) = \mathcal{G}\cdot (\mathcal{F}(x)) = \mathcal{G}\cdot (x\cdot\mathcal{F}) = \sum_i\sum_j\mathcal{G}_{i1} x_j\mathcal{F}_{ji}.
\end{split}
\end{align*}
Now, suppose that $\varphi\circ f=0$ for all $\varphi\in\Homf(V,\ef)$, i.e. $\varphi(f(x))= 0$ for any $x\in\ef^{(\beta)}$. Then
\begin{equation}\label{last}
\sum_i\sum_j\mathcal{G}_{i1}x_j\mathcal{F}_{ji} = 0
\end{equation}
for any $x\in\ef^{(\beta)}$ and any matrix $\mathcal{G}\in\ef^{\alpha,1}$. Given $j\in\beta$ and $i\in\alpha$, set for $\mathcal{G}$ the matrix such that
\begin{equation*}
\mathcal{G}_{ij} = \left\{
    \begin{array}{ll}
      0 & \hbox{if}\quad j\neq 1,\\
      1 & \hbox{if}\quad j =1.
    \end{array}
  \right.
\end{equation*}
Then \eqref{last} becomes
\begin{equation*}
\mathcal{F}_{ji} = 0
\end{equation*}
for all $(i,j)\in\alpha\times\beta$. Hence $f=0$.
\endproof

\end{document}